# Measurement of returns to scale in DEA using the CCR model

**Mahmood Mehdiloozad**

*Department of Mathematics, College of Sciences, Shiraz University, Shiraz 71454, Iran*
(E-mail: m.mehdiloozad@gmail.com, Tel.: +98.9127431689)

**Biresh K. Sahoo**

*Xavier Institute of Management, Xavier University, Bhubaneswar 751013, India*
(E-mail: biresh@ximb.ac.in, Tel.: + 91.6746647735)

(✉) **Corresponding author**: M. Mehdiloozad

Ph.D. Candidate

Department of Mathematics

College of Sciences

Shiraz University

Golestan Street | Adabiat Crossroad | Shiraz 71454 | Iran

# Measurement of returns to scale in DEA using the CCR model

**Abstract**

In data envelopment analysis (DEA) literature, the *returns to scale* (RTS) of an inefficient decision making unit (DMU) is determined at its projected point on the efficient frontier. Under the occurrences of multiple projection points, however, this evaluation procedure is not precise and may lead to erroneous inferences as to the RTS possibilities of DMUs. To circumvent this, the current communication first defines the RTS of an inefficient DMU at its projected point that lies in the relative interior of the minimum face. Based on this definition, it proposes an algorithm by extending the latest developed method of measuring RTS via the CCR model. The main advantage of our proposed algorithm lies in its computational efficiency.

## Introduction

In the data envelopment (DEA) literature on returns to scale (RTS), one of the approaches used for computing RTS is to look at the sum of the intensity variables at optimality of the CCR model of Charnes et al. (1978). The original form of this approach, as proposed by Banker and Thrall (1992), requires finding all optimal solutions of the CCR model, which is not practically feasible. Therefore, to address this crucial issue, Banker et al. (1996) proposed a two-stage method. For each decision making unit (DMU), the first stage applies the CCR model, to the DMU under evaluation if it is efficient or, to its projection point if it is inefficient. If the sum of the intensity variables at optimality is not equal to one, the second stage then obtains the lower or upper bound of this sum over all optimal solutions of the CCR model. In a relatively recent study, Zarepisheh et al. (2006) have enhanced the computational efficiency of Banker et al.'s (1996) method. Precisely, they have demonstrated that the RTS measurement of DMUs could be made by executing only the first stage of their method.

Note that since the RTS of an inefficient DMU is defined based on its projection point, it may not be determined *uniquely* when multiple projection points occur. Consequently, the definition of RTS will not be unambiguous if the projection point is chosen arbitrarily. Without consideration of this fact, any RTS measurement method may yield conflicting inferences on RTS possibilities for inefficient DMUs with multiple projections. For this reason, we believe that further improvement needs to be made to such methods.

Krivonozhko et al. (2012) have recently shown that all relative interior points of the *minimum face*—a face of minimum dimension that contains all the projection points—operate under the same type of RTS. This interesting finding reveals that the definition of RTS for an inefficient DMU can be made unambiguous by requiring its projection point to be in the relative interior of the associated minimum face. Based on this fact, we design a two-stage RTS measurement scheme by extending the algorithm of Zarepisheh et al. (2006). The first stage of our scheme uses the BCC model (Banker et al., 1984) for each DMU to identify both its efficiency status and projection point. The second stage measures the RTS at the DMU under evaluation via the envelopment form of the CCR model. If the DMU under evaluation is inefficient with a zero sum of optimal slacks, it is replaced by its BCC-projection point before proceeding to the second stage. However, if it is inefficient with a non-zero sum of optimal slacks, it is replaced by the projection point obtained from the linear programming model by Mehdiloozad et al. (2015a).

The main advantage of our proposed algorithm lies in its enhanced computational efficiency, which is due to two reasons. First, to find out the projection point of an inefficient DMU, the approach by Mehdiloozad et al. (2015a) is computationally more efficient than its alternatives. Second, to determine RTS, only a single DEA model is required to be solved in the envelopment form.

The rest of this communication is organized as follows. Section 2 contains some preliminaries, and Section 3 proposes our algorithm. Section 4 illustrates our proposed algorithm with a numerical example. Finally, Section 5 concludes with some remarks.

## Preliminaries

### Measurement of RTS via the CCR model

We deal with a technology set comprising of $n$ observed decision making units (DMUs), where each DMU$_j$ ($j \in J = \{1,...,n\}$) produces the output vector $\mathbf{y}_j = (y_{1j},...,y_{sj})^T \in \mathbb{R}_+^s$ by using the input vector $\mathbf{x}_j = (x_{1j},...,x_{mj})^T \in \mathbb{R}_+^m$. We denote the input and output matrices by $\mathbf{X} = [\mathbf{x}_1 \ ... \ \mathbf{x}_n]$ and $\mathbf{Y} = [\mathbf{y}_1 \ ... \ \mathbf{y}_n]$, respectively. We also use $o \in J$ as the index of the DMU under evaluation.

The CCR model and the BCC model for DMU$_o$ are set up, respectively, as:

$$\begin{aligned}
\min \quad & \theta - \varepsilon \left( \mathbf{1}_m^T \mathbf{s}^- + \mathbf{1}_s^T \mathbf{s}^+ \right) \\
\text{subject to} \quad & \\
& \mathbf{X}\boldsymbol{\lambda} + \mathbf{s}^- = \theta \mathbf{x}_o, \\
& \mathbf{Y}\boldsymbol{\lambda} - \mathbf{s}^+ = \mathbf{y}_o, \\
& \boldsymbol{\lambda} \geq \mathbf{0}_n, \ \mathbf{s}^- \geq \mathbf{0}, \ \mathbf{s}^+ \geq \mathbf{0},
\end{aligned} \quad (1)$$

and

$$\begin{aligned}
\min \quad & \theta - \varepsilon \left( \mathbf{1}_m^T \mathbf{s}^- + \mathbf{1}_s^T \mathbf{s}^+ \right) \\
\text{subject to} \quad & \\
& \mathbf{X}\boldsymbol{\lambda} + \mathbf{s}^- = \theta \mathbf{x}_o, \\
& \mathbf{Y}\boldsymbol{\lambda} - \mathbf{s}^+ = \mathbf{y}_o, \\
& \mathbf{1}_n^T \boldsymbol{\lambda} = 1, \\
& \boldsymbol{\lambda} \geq \mathbf{0}_n, \ \mathbf{s}^- \geq \mathbf{0}, \ \mathbf{s}^+ \geq \mathbf{0},
\end{aligned} \quad (2)$$

where $\varepsilon > 0$ is a positive non-Archimedean infinitesimal; and, $\mathbf{s}^-$ and $\mathbf{s}^+$ are the input and output slack vectors, respectively.

Let $\left(\theta^*, \boldsymbol{\lambda}^*, \mathbf{s}^{-*}, \mathbf{s}^{+*}\right)$ be an optimal solution to the CCR (BCC) model. Then, the optimal value $\theta^*$ is denoted by $\theta_o^{\text{CCR}}$ ($\theta_o^{\text{BCC}}$) and is called CCR- (BCC-) efficiency score.

Moreover, DMU$_o$ is said to be *CCR- (BCC-) efficient* if and only if $\theta_o^{CCR} = 1$ ($\theta_o^{BCC} = 1$), $\mathbf{s}^{-*} = \mathbf{0}_m$ and $\mathbf{s}^{+*} = \mathbf{0}_s$; otherwise, it is called BCC-inefficient.

Using the CCR model, Zarepisheh et al. (2006) presented the following theorem to measure RTS. Note that the evaluated DMU is assumed to be BCC-efficient; otherwise, it is replaced by its BCC-projection point defined by

$$\left(\mathbf{x}_o^{BCC}, \mathbf{y}_o^{BCC}\right) = \left(\theta_o^{BCC}\mathbf{x}_o - \mathbf{s}^{-*}, \mathbf{y}_o + \mathbf{s}^{+*}\right). \tag{3}$$

**Theorem 1** Let DMU$_o$ be a BCC-efficient unit. Then,

(i) constant RTS prevail at DMU$_o$ if and only if $\theta_o^{CCR} = 1$.

(ii) decreasing RTS prevail at DMU$_o$ if and only if $\theta_o^{CCR} < 1$ and $\mathbf{1}_n^T \boldsymbol{\lambda}^* > 1$ in any optimal solution of model (1).

(iii) increasing RTS prevail at DMU$_o$ if and only if $\theta_o^{CCR} < 1$ and $\mathbf{1}_n^T \boldsymbol{\lambda}^* < 1$ in any optimal solution of model (1).

Then, based on Theorem 1, they proposed the following algorithm to determine the RTS of DMU$_o$:

**Step 1:** Solve the BCC model. If DMU$_o$ is BCC-efficient, go to Step 2. Otherwise, first $(\mathbf{x}_o, \mathbf{y}_o) \leftarrow \left(\mathbf{x}_o^{BCC}, \mathbf{y}_o^{BCC}\right)$ and then go to Step 2.

**Step 2:** Solve the CCR model. If $\theta_o^{CCR} = 1$, then constant RTS prevail. Otherwise, decreasing RTS prevail if $\mathbf{1}_n^T \boldsymbol{\lambda}^* > 1$ and increasing RTS prevail if $\mathbf{1}_n^T \boldsymbol{\lambda}^* < 1$.

*Identification of the global reference set*

Let $\mathbf{X}_E$ and $\mathbf{Y}_E$ be respectively the input and output matrices of the BCC-efficient DMUs. Further, let $e$ be the cardinality of the set of all the BCC-efficient DMUs. We define $\Lambda_o$ as the set of all the BCC-projection points of DMU$_o$. Then, $\Lambda_o$ can be expressed as

$$\Lambda_o = \left\{ (\mathbf{X}_E\boldsymbol{\mu}, \mathbf{Y}_E\boldsymbol{\mu}) \middle| \begin{bmatrix} \mathbf{X}_E & \mathbf{I}_m & \mathbf{0}_{m\times s} \\ \mathbf{Y}_E & \mathbf{0}_{s\times m} & -\mathbf{I}_s \\ \mathbf{1}_e^T & \mathbf{0}_m^T & \mathbf{0}_s^T \\ \mathbf{0}_e^T & \mathbf{1}_m^T & \mathbf{1}_s^T \end{bmatrix} \begin{bmatrix} \boldsymbol{\mu} \\ \mathbf{s}^- \\ \mathbf{s}^+ \end{bmatrix} = \begin{bmatrix} \theta_o^{BCC}\mathbf{x}_o \\ \mathbf{y}_o \\ 1 \\ \mathbf{1}_m^T\mathbf{s}^{-*} + \mathbf{1}_s^T\mathbf{s}^{+*} \end{bmatrix}, \begin{bmatrix} \boldsymbol{\mu} \\ \mathbf{s}^- \\ \mathbf{s}^+ \end{bmatrix} \geq \mathbf{0}_{e+m+s} \right\}, \quad (4)$$

where $(\mathbf{s}^{-*}, \mathbf{s}^{+*})$ is a partial optimal solution to the BCC model.

We also define $\Omega_o$ as the set of all the intensity vectors that are associated with the BCC-projection points, i.e.,

$$\Omega_o := \left\{ \boldsymbol{\mu} \middle| (\mathbf{X}_E\boldsymbol{\mu}, \mathbf{Y}_E\boldsymbol{\mu}) \in \Lambda_o \right\}. \quad (5)$$

As has been demonstrated in Mehdiloozad et al. (2015a), the global reference set of $\mathrm{DMU}_o$, $R_o^G$, is a unique reference set containing all the possible reference units of $\mathrm{DMU}_o$; and can be identified from the following relation:

$$R_o^G = \left\{ (\mathbf{x}_j, \mathbf{y}_j) \middle| \mu_j^{\max} > 0 \right\}, \quad (6)$$

where $\boldsymbol{\mu}^{\max}$ is a *maximal element* of $\Omega_o$—an element with the maximum number of positive components.

Following Mehdiloozad et al. (2015a), $\boldsymbol{\mu}^{\max}$ can be found from the following linear programming model[1]:

$$\max \quad \mathbf{1}^T\boldsymbol{\mu}^1 + \delta^1$$
subject to
$$\begin{bmatrix} \mathbf{X}_E & \mathbf{I}_m & \mathbf{0}_{m\times s} \\ \mathbf{Y}_E & \mathbf{0}_{s\times m} & -\mathbf{I}_s \\ \mathbf{1}_e^T & \mathbf{0}_m^T & \mathbf{0}_s^T \\ \mathbf{0}_e^T & \mathbf{1}_m^T & \mathbf{1}_s^T \end{bmatrix} \begin{bmatrix} \boldsymbol{\mu}^1 + \boldsymbol{\mu}^2 \\ \mathbf{s}^- \\ \mathbf{s}^+ \end{bmatrix} - \begin{bmatrix} \theta_o^{BCC}\mathbf{x}_o \\ \mathbf{y}_o \\ 1 \\ \mathbf{1}_m^T\mathbf{s}^{-*} + \mathbf{1}_s^T\mathbf{s}^{+*} \end{bmatrix} (\delta^1 + \delta^2) = \mathbf{0}, \quad (7)$$
$$\mathbf{0}_e \leq \boldsymbol{\mu}^1 \leq \mathbf{1}_e, \quad 0 \leq \delta^1 \leq 1,$$
$$\boldsymbol{\mu}^2 \geq \mathbf{0}_e, \quad \mathbf{s}^- \geq \mathbf{0}_m, \quad \mathbf{s}^+ \geq \mathbf{0}_s, \quad \delta^2 \geq 0.$$

Precisely,

---

[1] Mehdiloozad (2015) demonstrated that this model was originally derived from a mixed 0-1 linear programming model by using the linear programming relaxation method.

$$\boldsymbol{\mu}^{\max} = \frac{1}{1+\delta^{2*}}\left(\boldsymbol{\mu}^{1*} + \boldsymbol{\mu}^{2*}\right), \tag{8}$$

where $\left(\boldsymbol{\mu}^{1*} + \boldsymbol{\mu}^{2*}, \mathbf{s}^{-*}, \mathbf{s}^{+*}, \delta^{1*}, \delta^{2*}\right)$ is an optimal solution to model (7).

## Our proposed algorithm

The observed BCC-inefficient DMUs can be divided into two groups. The first group ($G_1$) consists of the BCC-inefficient DMUs with zero sums of optimal slacks. The second group ($G_2$) consists of the BCC-inefficient DMUs with non-zero sums of optimal slacks for which multiple BCC-projection points may occur.

The standard approach followed in the DEA literature for measuring the RTS of an inefficient DMU is first to project it onto the BCC-efficiency frontier and then to determine its RTS at the BCC-projection point $\left(\mathbf{x}_o^{BCC}, \mathbf{y}_o^{BCC}\right)$ as defined in (3). If $G_2$ is empty, this approach is perfectly fine in defining RTS uniquely. Otherwise, this uniqueness property cannot be guaranteed since multiple BCC-projection points may reveal different types of RTS for the DMU under evaluation. Since the algorithm of Zarepisheh et al. (2006) is designed based on the standard approach, it may yield erroneous inferences on RTS possibilities for the units in $G_2$.

Therefore, to resolve this issue, the RTS must be defined over a subset of $\Lambda_o$ that its elements all exhibit the same RTS possibility. To accomplish the task, we resort to the concept of *minimum face*. As demonstrated by Krivonozhko et al. (2014), there exists a face of minimum dimension, called the minimum face, which contains $\Lambda_o$. On the other hand, Krivonozhko et al. (2012) have shown that all relative interior points of the minimum face operate under the same type of RTS. Thus, following Krivonozhko et al. (2014) and Mehdiloozad et al. (2015a), the RTS of a BCC-inefficient DMU is well defined over the intersection of $\Lambda_o$ with the relative interior of the minimum face. Based on this point, we present the following precise definition of RTS for the BCC-inefficient DMUs.

**Definition 1** The RTS of a BCC-inefficient DMU is defined at its projection point in the relative interior of its associated minimum face.

Note that a projection in the relative interior of the minimum face can be obtained by the following strict convex combination of the units in the global reference set:

$$\left(\mathbf{x}_o^{\max}, \mathbf{y}_o^{\max}\right) = \left(\mathbf{X}_E \boldsymbol{\mu}^{\max}, \mathbf{Y}_E \boldsymbol{\mu}^{\max}\right). \tag{9}$$

To determine the RTS possibility of DMU$_o$, we now propose the following two-stage scheme based on Definition 1:

**Stage 1:** Solve the BCC model for $(\mathbf{x}_o, \mathbf{y}_o)$.

- If $(\mathbf{x}_o, \mathbf{y}_o)$ is BCC-efficient, go to Step 2.
- Else if $(\mathbf{x}_o, \mathbf{y}_o) \in G_1$, obtain $\left(\mathbf{x}_o^{BCC}, \mathbf{y}_o^{BCC}\right)$ from (3). Then, replace $(\mathbf{x}_o, \mathbf{y}_o)$ with $\left(\mathbf{x}_o^{BCC}, \mathbf{y}_o^{BCC}\right)$ and go to Step 2.
- Else, solve model (7) and obtain $\left(\mathbf{x}_o^{\max}, \mathbf{y}_o^{\max}\right)$ from (9). Then, replace $(\mathbf{x}_o, \mathbf{y}_o)$ with $\left(\mathbf{x}_o^{\max}, \mathbf{y}_o^{\max}\right)$ and go to Step 2.

**Stage 2:** Solve the CCR model for $(\mathbf{x}_o, \mathbf{y}_o)$.

- If $\theta_o^{CCR} = 1$, then constant RTS prevail.
- Else if $\mathbf{1}_n^T \boldsymbol{\lambda}^* > 1$, decreasing RTS prevail.
- Else, increasing RTS prevail.

In summary, we estimate RTS of the observed DMUs by executing the following four main steps:

**Step 1** We first evaluate each DMU via the BCC model to obtain its efficiency score and the sum of the optimal input and output slacks. Based on the results, we then divide the BCC-inefficient DMUs into groups $G_1$ and $G_2$.

**Step 2** We determine RTS of the BCC-efficient DMUs.

**Step 3** We determine RTS of the DMUs in group $G_1$.

**Step 4** We determine RTS of the DMUs in group $G_2$.

**Remark 1** Note that after determining RTS of all the BCC-efficient DMUs in Step 2, the method of Tone (1996, 2005) can be applied in Step 4. This is worthwhile from the computational perspective because it avoids solving the CCR model.

**Remark 2** The most productive scale size (MPSS) pattern, as introduced by Banker (1984), determines the quantity by which a DMU with increasing (decreasing) RTS is expanded (contracted) in order to reach an optimal size. Davoodi et al. (2014) recently introduced the notion of *nearest MPSS* and defined it as "the closer MPSS pattern to the unit, the easier it would be to reach with that pattern" (p. 165). From their study, the MPSS pattern obtained from solving of the CCR model in Stage 2 of our proposed scheme is not guaranteed to be the nearest pattern. Moreover, to identify correctly the nearest MPSS pattern as well the RTS status of $(\mathbf{x}_o, \mathbf{y}_o)$, the following Stage 3 is added to our proposed scheme:

**Stage 3:** If decreasing (increasing) RTS prevail at $(\mathbf{x}_o, \mathbf{y}_o)$, then set $\frac{1}{\mathbf{1}_n^T \boldsymbol{\lambda}^*} (\theta_o^{CCR} \mathbf{x}_o, \mathbf{y}_o)$ as the nearest MPSS pattern, where $\boldsymbol{\lambda}^*$ is obtained by solving the following model:

$$\begin{aligned}
\min(\max) \quad & \mathbf{1}_n^T \boldsymbol{\lambda} \\
\text{subject to} \quad & \\
& \mathbf{X}\boldsymbol{\lambda} \leq \theta_o^{CCR} \mathbf{x}_o, \\
& \mathbf{Y}\boldsymbol{\lambda} \geq \mathbf{y}_o, \\
& \mathbf{1}_n^T \boldsymbol{\lambda} \geq 1 \; (\mathbf{1}_n^T \boldsymbol{\lambda} \leq 1), \\
& \boldsymbol{\lambda} \geq \mathbf{0}_n,
\end{aligned} \quad (10)$$

**A numerical example**

Let us consider a two-inputs–one-output technology characterized by six hypothetical DMUs labeled as A–F. The observed input–output data are all exhibited in Table 1 and the resulting BCC technology set is depicted in Fig. 1. With this example, we first show how the occurrence of multiple projection points yields different RTS possibilities for an inefficient DMU; and then show how our proposed scheme deals effectively with this issue.

**Table 1** Input-output data for six DMUs

| DMU | $x_1$ | $x_2$ | $y$ |
|-----|-------|-------|-----|
| A | 4 | 1 | 1 |
| B | 4 | 2 | 2 |
| C | 6 | 1 | 3 |
| D | 9 | 1.5 | 3 |
| E | 4 | 2 | 1 |
| F | 9 | 1.5 | 1 |

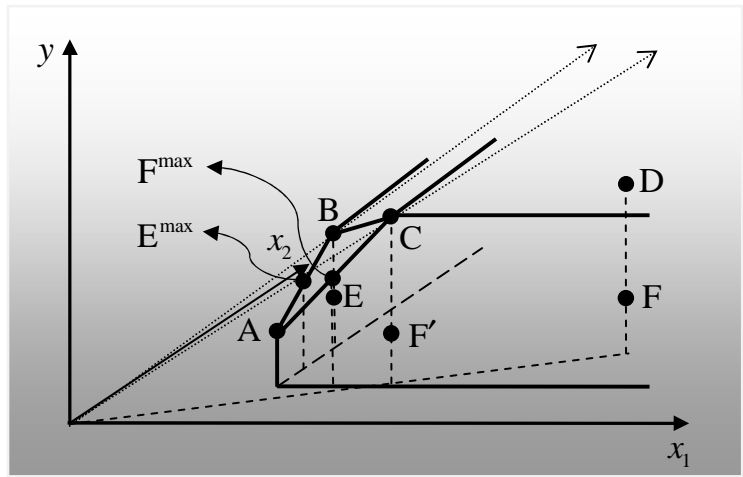

**Fig. 1** The technology set spanned by units A–F

***Step 1:*** *Evaluation of DMUs via the BCC model*

We first evaluated all the units via the BCC model and then presented the results in Table 2.

For each DMU, Table 2 exhibits the BCC-efficiency score, the optimal sum of slacks and the corresponding BCC-projection point. The last column of this table also indicates the group to which each BCC-inefficient DMU is assigned. As can be seen from Table 2 and Fig. 1, units A, B and C are all BCC efficient, and units D, E, and F are BCC inefficient.

**Table 2** The results obtained from the BCC model

| DMU | $\theta_o^{BCC}$ | $\mathbf{1}_m^T \mathbf{s}^{-*} + \mathbf{1}_s^T \mathbf{s}^{+*}$ | $x_1^{BCC}$ | $x_2^{BCC}$ | $y^{BCC}$ | group |
|---|---|---|---|---|---|---|
| A | 1 | 0 | -- | -- | -- | -- |
| B | 1 | 0 | -- | -- | -- | -- |
| C | 1 | 0 | -- | -- | -- | -- |
| D | 0.6667 | 0 | 6 | 1 | 3 | $G_2$ |
| E | 1 | 1 | 4 | 2 | 2 | $G_2$ |
| F | 0.6667 | 2 | 4 | 1 | 1 | $G_1$ |

As can be seen in Fig. 1, the BCC model first radially contracts the inputs of units D and F. While unit D is projected onto the BCC-efficient unit C, unit F is projected onto the frontier at the BCC-inefficient point $F' = ((6,1)^T, 1)$. By decreasing the first input of $F'$ by two units, the BCC model then determines unit A as the BCC-projection point of unit F. Note that the points on the line segment AC are all the possible projection points for unit F. For example, an increase in the output of $F'$ by two units results in unit C as another projection point of unit F. Moreover, the minimum face is the line segment AC itself. Similarly, the projection set and the minimum face for unit E are the same, i.e., the line segment AB, and its BCC-projection point is unit B. In summary, the results reveal that $G_1 = \{F\}$ and $G_2 = \{D, E\}$.

*Step 2: RTS Measurement for the BCC-efficient units*

From the results of Step 1, we found units A, B, and C as BCC efficient. Using the CCR model, we evaluated these units and reported the results in Table 3. Since the CCR efficiency scores for units B and C are all equal to one, these units operate under constant RTS. For better exposition, we depicted, in Fig. 1, the rays passing through these units.

Since for unit A, we have $\theta_A^{CCR} = 0.5 < 1$ and $\mathbf{1}_n^T \boldsymbol{\lambda}^* = 0.3750 < 1$, Theorem 1 implies that this unit exhibits increasing RTS. Hence, following Tone (1996, 2005), all activities on the triangular ABC, excepting those on the line segment BC, exhibit increasing RTS. This is because the global reference set of all such activities consist of units A, B, and C.

**Table 3** The RTS of the BCC-efficient units

| DMU | $\theta_o^{CCR}$ | $\mathbf{1}_n^T \boldsymbol{\lambda}^*$ | RTS |
|---|---|---|---|
| A | 0.5 | 0.3570 | I |
| B | 1 | -- | C |
| C | 1 | -- | C |

**Note:** C: Constant RTS, and I: Increasing RTS

*Step 3: RTS Measurement for the units in $G_1$*

In the first stage, unit F is found to be the only unit in $G_1$. Since the BCC-projection point of unit F is unit C, this DMU exhibits constant RTS.

*Step 4: RTS Measurement for the units in $G_2$*

We now turn to determine the RTS possibilities of units E and F, which were all identified in Step 1 as the BCC-inefficient units with non-zero optimal sums of slacks. To accomplish this, we solved model (7) for these two units and reported the result in Table 4. For each DMU, Table 4 exhibits the reference units together with their associated weights, the projection points in the relative interior of the minimum face and the RTS status.

**Table 4** The RTS of the units of group $G_2$

| | Global reference set | | | projection | | | |
|---|---|---|---|---|---|---|---|
| DMU | A | B | C | $x_1^{max}$ | $x_2^{max}$ | $y^{max}$ | RTS |
| E | 0.5 | 0.5 | -- | 4 | 1.5 | 1.5 | I |
| F | 0.5 | -- | 0.5 | 5 | 1 | 2 | I |

**Note:** Increasing RTS

Let us first consider unit E. As can be seen from Table 4, the reference units for this DMU are units A and B with weights 0.5 and 0.5, respectively, thereby yielding the projection point $E^{max} = ((4,1.5)^T, 1.5)$. This projection lies obviously in the relative interior of the associated minimum face, i.e., the line segment AB, and exhibits increasing RTS.

Hence, based on Definition 1, the precise RTS status of unit E is increasing. However, the RTS status of unit E based on its BCC-projection (i.e., unit B) is constant.

Now, we consider unit F whose reference units are A and C, each with the respective weight of 0.5. Using model (7) we obtained its projection as $D^{max} = \left( (5,1)^T, 2 \right)$, which lies in the relative interior of the line segment AC, thus exhibiting increasing RTS. Hence, as per our Definition 1, the precise RTS status of unit F is increasing. This is *accidentally* the same RTS found based on the BCC-projection (i.e., unit A).

## Concluding remarks

The set of all the BCC-inefficient DMUs can be divided into two disjoint groups $G_1$ and $G_2$, which, respectively, contain DMUs with zero and non-zero sums of optimal slacks. While each unit in $G_1$ is uniquely projected onto the BCC-efficiency frontier, multiple projection points may occur for each unit of group $G_2$. This indicates that the concept of RTS has an ambiguous meaning for the DMUs of group $G_2$. This is because the RTS of a BCC-inefficient unit is determined at its projection point and, consequently, multiple projection points may reveal different RTS possibilities for this DMU. Therefore, this ambiguity may lead to erroneous inferences as to the RTS possibilities of DMUs and thereby, it adversely influences the identification of the nearest MPSS pattern.

In this communication, the above-mentioned ambiguity is effectively eliminated by defining the RTS over a subset of the projection set that its elements all exhibit the same RTS possibility. Precisely, the RTS of a BCC-inefficient DMU is defined over the intersection of the projection set with the relative interior of the minimum face. This definition is precise in accordance with the interesting finding of Krivonozhko et al. (2012) that all relative interior points of the minimum face exhibit the same RTS possibility.

Using this new definition of RTS, a two-stage scheme is then proposed for determining the RTS possibilities of DMUs. For the BCC-efficient and the BCC-inefficient units in $G_1$, the operations of the proposed scheme and those of Zarepisheh et al.'s (2006) algorithm are the same. That is, the proposed scheme applies the CCR model, to the DMU under evaluation if it is efficient or, to its projection point if it is inefficient. For the BCC-inefficient units of group $G_2$, however, the proposed scheme applies the linear programming model of Mehdiloozad et al. (2015) to the DMU under evaluation to obtain a projection point that lies in the relative interior of the minimum face. It then determines the RTS at this projection point. Toward precise identification of the nearest MPSS pattern as well the RTS status of the DMU under evaluation, a third stage is also added to the proposed scheme.

Note that, to identify the projection point of a DMU in group $G_2$, the approach by Mehdiloozad et al. (2015a) is computationally more efficient than its alternatives. On the other hand, to determine RTS, only a single DEA model is required to be solved in the envelopment form. These two points demonstrate the enhanced computational efficiency of our proposed scheme.

Finally, due to the importance of determining RTS in weight-restricted DEA models, the extension of our proposed scheme to the weight-restricted DEA framework is suggested as a future research subject—in the spirit of Hosseinzadeh Lotfi et al. (2007) and Mehdiloozad et al. (2015b).